\documentstyle[12pt]{article}

  \newcommand{\const}{\rm const}
  \newcommand{\Var}{\rm Var}

\begin{document}

 \begin{center}

 \ {\bf Signed variable  optimal kernel for} \par

 \vspace{4mm}

 {\bf non - parametric density estimation. }\par

\vspace{5mm}

{\bf  M.R.Formica, E.Ostrovsky, and L.Sirota. } \par

 \vspace{5mm}

 \end{center}

 \ Universit\`{a} degli Studi di Napoli Parthenope, via Generale Parisi 13, Palazzo Pacanowsky, 80132,
Napoli, Italy. \\

e-mail: mara.formica@uniparthenope.it \\

\vspace{3mm}

Department of Mathematics and Statistics, Bar-Ilan University,
59200, Ramat Gan, Israel. \\

\vspace{4mm}

e-mail: eugostrovsky@list.ru\\
Department of Mathematics and Statistics, Bar-Ilan University,\\
59200, Ramat Gan, Israel.

\vspace{4mm}

e-mail: sirota3@bezeqint.net \\

\vspace{5mm}

\begin{center}

 \ {\bf Abstract.}

 \end{center}

  \ We derive the optimal signed variable in general case kernels for the classical statistic density estimation, which are some
 generalization of the famous  Epanechnikov's  ones.\par

 \vspace{5mm}

  \par

\vspace{5mm}

\section{Statement of problem. Notations and definitions. Previous results.}

\vspace{5mm}

 \ Let  $ \ (\Omega,  M, {\bf P})  \ $  be probability space with expectation $ \ {\bf E} \ $ and variance $ \ \Var. \ $
 \ Let also  $ \ \{\xi_k \}, \ k = 1,2,\ldots,n \ $ be a sequence of independent, identical distributed (i, i.d.) random variables (r.v.)
 taking the values in the  real axis $ \ R, \ $ and having certain {\it non - known} density of a distribution $ \ f = f(x), \ x \in R. \ $ \par
  \ We suppose further that this function belongs to the certain space of all the numerical valued $ \  2m \ -  \ $ times continuous bounded differentiable
 functions  $ \ C^{2m}(R), \ m = 1,2,\ldots,   \ $ having a finite norm

\begin{equation}  \label{C2m space}
||f||_{2m} \stackrel{def}{=} \max_{k= 0,1,2,\ldots,2m} \ \sup_{x \in R} | f^{(k)}(x) | < \infty.
\end{equation}

\  Let also $ \ K = K(x), \ x \in  \ $ be certain {\it kernel},  i.e. measurable {\it even} function having finite support:

\begin{equation} \label{theta support}
\exists \theta \in (0,\infty) \ \forall x: |x| > \theta \Rightarrow K(x) = 0,
\end{equation}

 for which

\begin{equation} \label{Kernel}
\int_{R} K(x) \ dx = 1.
\end{equation}

 \ We impose also the following conditions on this kernel.

\begin{equation} \label{K cond 1}
K(-x) = K(x);  \ V_2(K) := \int_{R} K^2(x) dx < \infty;
\end{equation}

\begin{equation} \label{K cond 2}
K(\cdot) \in C(R), \ \int_{R} |K(x)| dx < \infty.
\end{equation}

 \ The following conditions, which are also to be presumed, may be named as  {\it conditions of orthogonality:}

\begin{equation} \label{first cond orthog}
\forall l = 1,2,\ldots,2m - 1 \ \Rightarrow \int_{R} x^l \ K(x) \ dx = 0.
\end{equation}

 \ Recall that the classical Parzen - Rosenblatt estimation $ \ f_n(x) = f_n^{PR}(x)  \ $
 of the density function $ \ f(x) \ $ has a form

\begin{equation} \label{OPR}
f_n(x) \stackrel{def}{=} \frac{1}{nh} \ \sum_{i=1}^n K \left( \ \frac{x - \xi_i}{h}  \ \right).
\end{equation}
 \ Here $ \ h = h(n)  \ $  be a deterministic positive sequence such that $ \ \lim_{n \to \infty} h(n) = 0 \ $
 and $\ \lim_{n \to \infty} n h(n) = \infty, \ $  see \cite{Parzen}, \cite{Rosenblatt}. \par

 \ These and alike estimations  was study in many works, see e.g. \cite{Devroe},  \cite{Epanechnikov},  \cite{Ibragimov Khasminskii 1},
 \cite{Ibragimov Khasminskii 2},  \cite{Nadaraya}, \cite{Parzen}, \cite{Rosenblatt},
 \cite{Slaoui},  \cite{Tsybakov1}, \cite{Tsybakov2}, \cite{Wegman}, \cite{Wolverton Wagner} etc. The optimal choose of
$ \ h = h(n)  \ $ and the kernel $ \ K(x) \ $ are devoted the following works
\cite{Compte},   \cite{Goldenshluger Lepski},  \cite{Khardani}, \cite{Lerasle}. The case when the r.v. - s. $ \ \{ \ \xi_i \ \} \ $
 are (weakly) dependent is investigated in  \cite{Khardani}, \cite{Simone}. \par

\  The conditions (\ref{first cond orthog})  may be used only for the investigation of {\it bias} $ \ \delta_n(x) \ $ of these statistics.
Namely, denote

$$
\delta_n(x) \stackrel{def}{=} {\bf E} f_n(x) - f(x);
$$
then under our conditions

\begin{equation} \label{bias est}
|\delta_n(x)| \le C_1(f) h^{2m}(n).
\end{equation}

 \ The variation of $ \ f_n(x) \ $ may be estimated as follows

\begin{equation} \label{Var es}
\Var \{ \ f_n(x) \ \} \le C_2(f) \ n^{-1} \ \int_{-\infty}^{\infty} K^2 (y) \ dy.
\end{equation}

\vspace{4mm}

\begin{center}

{\sc Statement of an optimization problem.}

\end{center}

\vspace{4mm}

 \ The relations   (\ref{Var es}) common with the limitations  (\ref{theta support}),   (\ref{Kernel}), (\ref{K cond 1}), (\ref{K cond 2}),
(\ref{first cond orthog}) lead as was shown  by V.A.Epanechnikov  in \cite{Epanechnikov} to the following setting of the constrained  optimization problem:

\begin{equation} \label{statement}
\int_{-\infty}^{\infty} K^2(y) \ dy \to \min,
\end{equation}

\vspace{3mm}

under limitations

\vspace{3mm}

\begin{equation} \label{restr1}
\int_{-\infty}^{\infty} K(y) \ dy = 1; \ \forall l=1,2,3, \ldots,2m - 1 \ \Rightarrow \int_{-\infty}^{\infty} y^l \ K(y) \ dy = 0;
\end{equation}

\begin{equation} \label{2m moment}
\int_{-\infty}^{\infty} y^{2m} K(y) dy = 1,
\end{equation}

\begin{equation} \label{theta restr}
\exists \theta \in (0,\infty) \  \forall y: |y| > \theta \  \Rightarrow K(y) = 0.
\end{equation}

 \vspace{5mm}

 \ This problem was solved in the case $ \ m = 1 \ $ by V.A.Epanechnikov in \cite{Epanechnikov}, 1969 year: $ \ \theta = \sqrt{5} \ $ and

\begin{equation} \label{Epanechnikov solution}
|y| \le \sqrt{5} \Rightarrow K(y) = \frac{3}{4 \sqrt{5}} - \frac{3y^2}{20 \sqrt{5}}.
\end{equation}

\vspace{4mm}

 \  {\bf  Our aim in this short report is to find the optimal kernel $ \ K \ $ for arbitrary natural value } $ \ m = 2,3,\ldots.\ $ \par
{\bf The case when the number $ \ m \ $ is fractional, will be considered the fourth section.} \par

\vspace{5mm}

 \ Another motivation of these statement of problem appears from the famous result belonging to W.Stute \cite{Stute 1}, \cite{Stute 2}:

$$
\lim_{n \to \infty} \sqrt{ \ \frac{n h(n)}{2 |\ln h(n)|} \ }  ||f_n - {\bf E} f_n||_{\infty} = ||f||_{\infty}^{1/2} \cdot V_2(K),
$$
where as ordinary $ \ ||g||_{\infty} = \sup_{x \in R} |g(x)|. \ $\par

\vspace{4mm}

 \ {\bf Remark 1.1.} When $ \ m \ge 2, \ $ we imposed in particular the following condition on the kernel $ \ K(\cdot): \ $

$$
\int_{-\theta}^{\theta} y^2 \ K(y) \ dy = 0.
$$
 \ Therefore, the kernel $ \ K(\cdot) \ $ can not take only non - negative values. \par

\vspace{5mm}

\section{Main result.}

\vspace{5mm}

\begin{center}

 {\sc Some facts  about Legendre's polynomials.    }

\end{center}

  \ The classical Legendre's  polynomials with support on the closed interval $ \ X := [-1,1], \ $ denotes as ordinary by
$ \ P_k(x), \ x \in [-1,1], \ k = 0,1,2,\ldots  \ $ may be defined for instance as follows

$$
P_k(x) = \frac{1}{2^k \ k!} \ \frac{d^k}{d x^k} (x^2 - 1)^k,
$$
Rodriguez'formula. The polynomial $ \ P_k(x) \ $ is really polynomial of  degree $ \ k. \ $   These polynomials are orthogonal:

$$
\int_{-1}^1 P_k(x) \ P_l(x) \ dx  = \frac{2}{2k +1} \ \delta_{k,l}, \ k,l = 0,1,2,\ldots,
$$
where $ \ \delta_{k,l} \ $ is Kroneker's symbol. Following,

\begin{equation} \label{orthogon}
\forall k \ge 1, l < k, \ l \ge 0 \ \Rightarrow \int_{-1}^1 P_k(x) \ x^l \ dx = 0.
\end{equation}
\ Many properties of these polynomials may be found, e.g.  in  the classical book \cite{Bateman}. We will use
the following relations:

\begin{equation} \label{value at one}
\forall k = 0,1,2,\ldots \ P_k(1) = 1;
\end{equation}

\begin{equation} \label{import int}
 2\mu(k) \stackrel{def}{=} \int_{-1}^1 x^k \ P_k(x) \ dx = \frac{2^{k+1} \ (k!)^2}{(2k + 1)!}.
\end{equation}

 \ Some examples: $ \ P_0(x) = 1,\ P_1(x) = x, \ P_2(x) = 0.5 (3x^2 - 1), \ $

\begin{equation} \label{P4}
 P_4(x) = 8^{-1} (35  x^4 - 30 x^2 +3).
\end{equation}

\vspace{4mm}

 \ {\bf Definition 2.1.} Let $ \ \theta = \const \in (0,\infty). \ $ The {\it dilated} Legendre's polynomial
 $ \ L_k^{\theta}(x), \ k=0,1,2,\ldots; \ x \in [-\theta,\theta] \ $ of degree $ \ k, \ k = 0,1,2,\ldots \ $
 with parameter $ \ \theta \ $  is defined as follows.

\begin{equation} \label{dilated pol}
L_k^{\theta}(x) \stackrel{def}{=} \frac{1}{\theta} \ P_k \left( \ \frac{x}{\theta}  \ \right).
\end{equation}

\vspace{3mm}

 \ Of course, the properties of these polynomials  follows from ones for Legendre's polynomials. For instance,
$ \ L_{2s}(\cdot), s = 0,1,2,\ldots \ $ is even function,  $ \ L_{2s+1}(\cdot), s = 0,1,2,\ldots \ $ is odd;  \\
relations of orthogonality

\begin{equation} \label{ortho dilat}
\int_{-\theta}^{\theta} x^l \ L_k^{\theta}(x) \ dx = 0, \ k > l, \ l = 0,1,\ldots, k-1;
\end{equation}
 as well as

$$
 L_k^{0}(0) = 1/\theta = L_k(\theta),
$$

\begin{equation} \label{L import int}
\int_{-\theta}^{\theta} y^k \ L_k^{\theta}(y) \ dy =  \theta^k \cdot \frac{2^{k+1} \ (k!)^2}{(2k + 1)!},
\end{equation}

\begin{equation} \label{L orthogon}
\int_{-\theta}^{\theta} L_k(y) \ L_l(y) \ dy = \frac{2/\theta}{2k + 1} \ \delta_{k,l}.
\end{equation}
\ In particular,

\begin{equation} \label{L norming}
\int_{-\theta}^{\theta} L_k^2(y) \ dy = \frac{2/\theta}{2k + 1}.
\end{equation}

 \vspace{4mm}

  \ Let us return to the formulated above optimization problem  (\ref{statement}), (\ref{restr1}), (\ref{2m moment}), (\ref{theta restr}). \par

\vspace{4mm}

 \ {\bf Theorem 2.1.} The formulated optimization problem has an unique solution $ \ K_0(y) \ $ and having a
 following form

\begin{equation} \label{optimal 2m kernel}
K_0(y) = \frac{1}{2 \theta_0} - \frac{1}{2} L^{\theta_0}_{2m}(y) =  \frac{1}{2 \theta_0} -  \frac{1}{2 \theta_0} P_{2m} \left( \ \frac{y}{\theta_0}  \ \right),
\end{equation}

when $ \ |y| \le \theta_0, \ $ and $ \ K_0(y) = 0 \ $ otherwise. Here

\begin{equation} \label{theta 0}
\theta_0 = \left[ \  \frac{1}{1 -  \mu(2m)} \ \right]^{1/(2m)},
\end{equation}

\vspace{3mm}

and wherein under our restrictions on the kernel $ \ K(\cdot) \ $

\begin{equation} \label{minimal value}
\min \int_R K^2(y) dy =  \int_R K^2_0(y) \ dy  = \frac{1}{2 \theta_0} \cdot \frac{4m+3}{4m+1}.
\end{equation}

\vspace{3mm}

 \ Recall that

$$
\mu(2m) = \frac{2^{2m} \ (2m)!}{(4 m + 1)!}.
$$

\vspace{4mm}

 \ {\bf Proof.} We will follow the V.A.Epanechnikov,  author of the article \cite{Epanechnikov}, applied the famous
 Legendre factors method. The Euler's equation for this problem give us the following equality for the optimal kernel $ \ K(\cdot) \ $

 $$
 K(y) = \sum_{s=0}^{2m} \lambda_s \ y^s, \ |y| \le \theta.
 $$
 \ On the other words, $ \ K \ $ is polynomial of degree $ \ \le 2m, \ $  inside the interval $ \ [-\theta,\theta]. \ $
As long as the kernel is even function,

$$
\lambda_{2 r + 1} = 0, \ r = 0,1,\ldots, m - 1.
$$
 \ Further, it follows from the conditions of orthogonality that also

$$
\lambda_2 = \lambda_4 = \ldots \lambda_{2m - 2} = 0.
$$
 \ Therefore, the optimal kernel has a form

$$
K(y) = a - b L^{\theta}_{2m}(y) = a - \frac{b}{\theta} P_{2m} \left( \ \frac{y}{\theta} \ \right), \ a,b = \const.
$$
 \ We deduce substituting  the value $ \ y = \theta \ $ and taking into account the relations $ \ K(\theta) = 0 \ $
 and $ \ P_{2m}(1) = 1  \ $

\begin{equation} \label{first eq}
a = \frac{b}{\theta}.
\end{equation}

 \ Secondly,

 $$
 1 = \int_{-\theta}^{\theta} K(y) dy = 2 a \theta,
 $$
following $ \ a = 1/(2\theta) \ $ and hence $ \ b = 1/2, \ $ so that

$$
K(y) = \frac{1}{2 \theta} - \frac{1}{2} L_{2m}^{\theta} (y) = \frac{1}{2 \theta} - \frac{1}{2 \theta} P_{2m} \left( \ \frac{y}{\theta}  \ \right).
$$

\vspace{3mm}

 \ Thirdly,

$$
 1 = \int_{-\theta}^{\theta} y^{2m} \ K(y) \ dy  = \theta^{2m} - 0.5 \ \theta^{2m} \int_{-1}^1 z^{2m} \ P_{2m}(z) \ dz =
$$

$$
\theta^{2m} - \theta^{2m} \ \mu(2m).
$$

\vspace{3mm}

 \ Now, proposition (\ref{theta 0}) there holds, as well. Thus,

$$
K_0(y)= \frac{1}{2 \theta_0} -  \frac{1}{2 \theta_0} P_{2m} \left( \ \frac{y}{\theta_0}  \ \right), \ |y| \le \theta_0,
$$

 \ Ultimately, the equality (\ref{minimal value}) follows immediately from ones (\ref{L orthogon}) and (\ref{L norming}).\par

\vspace{3mm}

 \ This completes the proof of our theorem. \par

\vspace{5mm}

 \ {\bf Examples.} Of course, for the case $ \ m = 1 \ $ we obtain the classical Epanechnikov's kernel. \par
  \ Let now $ \ m = 2. \ $ We deduce after some calculations

$$
\theta_0 = \left[ \ {\frac{63}{11}} \ \right]^{1/4},
$$
and correspondingly

$$
K_0(y) = \frac{1}{2\theta_0} \ \left\{ \ 1 - \frac{1}{8} \left( \ 35 y^4 \theta_0^{-4} -30 y^2 \theta_0^{-2} +3 \ \right)  \ \right\}, \
|y| \le \theta_0,
$$
 and of course

$$
K_0(y) =0, \ |y| > \theta_0.
$$

 \ Wherein

\begin{equation} \label{minimal part  value}
\min \int_R K^2(y) \ dy =  \int_R K^2_0(y) \ dy  = \frac{11}{18 \ \theta_0}.
\end{equation}

 \vspace{5mm}

\section{The case of fractional order.}

 \vspace{5mm}

 \ Let $ \ \beta \ $ be arbitrary  {\it  fractional positive} number; denote by $ \ l = [\beta] \ $ its positive integer part.
We suppose that the function $ \ f(\cdot) \ $ belongs to the space $ \ \Sigma(\beta). \ $ This imply that all the derivatives
$ \ f^{(k)}(x), k = 0,1,\ldots,l \ $ are continuous and bounded and that the last continuous derivative $ \ f^{(l)}(x) \ $ is
bounded  and satisfies the H\"older's condition with power $ \ \beta - l. \ $ \par
 \ We impose on the kernel $ \ K(\cdot) \ $ conditions alike ones in the first section:  $ \ K = K(x), \ x \in  \ $ be certain
 even function having finite support:

\begin{equation} \label{theta beta support}
\exists \theta \in (0,\infty) \ \forall x: |x| > \theta \Rightarrow K(x) = 0,
\end{equation}

 for which

\begin{equation} \label{Kernel beta}
\int_{R} K(x) \ dx = 1,
\end{equation}

\begin{equation} \label{K cond  beta 1}
K(-x) = K(x); \ \int_{R} K^2(x) dx < \infty;
\end{equation}

\begin{equation} \label{K  cond beta 2}
K(\cdot) \in C(R), \ \int_{R} |K(x)| dx < \infty.
\end{equation}

 \ The following conditions may be named as  before {\it conditions of orthogonality:}

\begin{equation} \label{cond orthog}
\forall r = 1,2,\ldots, l \ \Rightarrow \int_{R} x^r \ K(x) \ dx = 0.
\end{equation}

 \ Ultimately,  suppose

\begin{equation} \label{beta mom}
\int_R |x|^{\beta} K(x) \ dx < \infty.
\end{equation}

 \ It is known that under these conditions the bias $ \ \delta_n(x) \ $ of Parzen - Rosenblatt's estimation  $ \ f_n(x), \ $  as well as
 for Wolverton - Wagner's  ones  obey's a following  property

$$
\sup_x |\delta_n(x)| = \sup_x | \ {\bf E} f_n(x) - f(x) \ | \le C(\beta,f) h_n^{\beta},
$$
see e.g. \cite{Compte},  \cite{Ibragimov Khasminskii 1}, \cite{Ibragimov Khasminskii 2}. \par
 \ We get following again V.A.Epanechnikov \cite{Epanechnikov} the following variational statement of problem
under formulated above in this section restrictions

\begin{equation} \label{beta optimiz}
\int_R K^2(y) \ dy \to \min_K
\end{equation}
where {\it in addition} $ \ \int_R |y|^{\beta} \ K(y) \ dy = 1.\ $ \par

\vspace{4mm}

 \ {\bf Theorem 3.1.} The (unique)  solution of this problem has a form

\begin{equation} \label{theta0 beta}
 \theta_0 = (2 \beta + 1)^{1/\beta};
\end{equation}

\begin{equation} \label{K0 beta}
K_0(y) = \lambda - \mu |y|^{\beta}, \ |y| \le \theta_0, \ K_0(y) = 0, \ |y| > \theta_0;
\end{equation}

\begin{equation} \label{lambda beta}
\lambda =  \frac{\beta + 1}{2\beta} \cdot (2 \beta + 1)^{-1/\beta},
\end{equation}

\begin{equation} \label{lambda beta}
\mu = \frac{\beta + 1}{2\beta} \cdot (2 \beta + 1)^{-(\beta + 1)/\beta},
\end{equation}

and wherein under our condition

\begin{equation} \label{minimal value 1}
\min \int_R K^2(y) \ dy = \int_R K_0^2(y) \ dy = \mu^2 \ \theta_0^{2 \beta + 1} \ \cdot \left\{ \ \frac{2 \beta^2}{(2 \beta + 1)(\beta +1)}  \ \right\}=
\end{equation}

\begin{equation} \label{minimal value 2}
(\beta +1) \cdot (2 \beta + 1)^{- (\beta + 1)/\beta}.
\end{equation}

\vspace{4mm}

 \ {\bf Proof.} The Euler's equations  for this problem give us as before the following form for the
optimal kernel

\begin{equation} \label{form of K}
K(y) = \lambda - \mu |y|^{\beta} + \sum_{j=1}^l  \nu_j y^j, \ |y| \le \theta.
\end{equation}
 \ Since the function $ \ K(\cdot) \ $ is even, $ \ \nu_1 = \nu_3 = \ldots  = 0.  \ $  Further, it follows from the
relations of orthogonality that all the other coefficients  $ \ \nu_j \ $ are absent; so the optimal kernel $ \ K \ $
has a form inside the closed interval $ \ y \in [ - \theta,\theta]: \  $

\begin{equation} \label{form of K}
K(y) = \lambda - \mu \ |y|^{\beta}.
\end{equation}

\ As long as $ \ K(\theta) = 0, \ $

\begin{equation} \label{first eq}
\lambda = \mu \cdot \theta^{\beta}.
\end{equation}

\vspace{3mm}

\ Secondly,

$$
1 = \int_{-\theta}^{\theta} K(y) \ dy =  2 \left[ \ \lambda \theta - \mu \frac{\theta^{\beta}}{\beta + 1}   \ \right],
$$

following

\begin{equation} \label{second eq}
\lambda \theta -   \mu \cdot \frac{\theta^{\beta}}{\beta + 1}  = \frac{1}{2}.
\end{equation}

\vspace{3mm}

 \ Thirdly,

$$
1 = \int_{-\theta}^{\theta} |y|^{\beta} K(y) \ dy =  2 \ \lambda \frac{\theta^{\beta}}{\beta + 1} - 2 \ \mu \frac{\theta^{2 \beta + 1}}{2 \beta + 1},
$$
 or equally

\begin{equation} \label{third eq}
\lambda \frac{\theta^{\beta + 1}}{\beta + 1} - \mu \frac{\theta^{2 \beta +1}}{2 \beta + 1}  = \frac{1}{2}.
\end{equation}

\vspace{3mm}

 \ Solving  the system of equations (\ref{first eq}), \  (\ref{second eq}) and  (\ref{third eq}), we obtain the assertion of theorem 3.1. \par

\vspace{3mm}

 \ The equalities  (\ref{minimal value 1}), (\ref{minimal value 2}) may be obtained after simple calculations.\par

\vspace{3mm}

 \ When for instance $ \ \beta = 3/2, \ $ then

$$
\theta_0 = 4^{2/3}, \ \lambda = \frac{5}{6} \cdot 4^{-2/3}, \ \mu =   \frac{5}{6} \cdot 4^{-5/3},
$$

and

$$
\min \int_R K_0^2(y) \ dy = 5 \cdot 2^{-1/3}.
$$

\vspace{5mm}

 \ {\bf Remark 3.0.} \par

  \ Note that in the case of fractional value $ \ \beta \ $ the optimal kernel $ \ K_0 \ $ is non - negative!

\vspace{5mm}

 \ {\bf Remark 3.1.} \par

 \ It is  interest to note that if we  choose as the value $ \ \beta \ $ an {\it integer} value $ \ \beta := 2, \ $ we obtain the classical
 Epanechnikov's kernel

$$
\theta_0 = \sqrt{5}, \  \lambda = \frac{3}{4 \sqrt{5}}, \ \mu = \frac{3}{20 \sqrt{5}}.
$$

 \ In this case the minimal value of $ \ \int_R K_0^2(y) \ dy \ $ is equal to $ \ 3/(5 \sqrt{5}). \ $

\vspace{4mm}

\ {\bf Remark 3.2.} We do not use in this section, i.e. in the case of fractional value of the parameter $ \ \beta, \ $
the theory of Legendre's polynomials, in contradiction to the foregoing sections. \par

\vspace{5mm}

\section{Another statement of problem.}

 \vspace{5mm}

 \ We retain all the restrictions and conditions of the foregoing section.  Denote  in addition

$$
J_{\beta} = J_{\beta}(K) := \int_R |y|^{\beta} \ K(y) \ dy; \hspace{4mm}  V_2 = V_2(K) := \int_R K^2(y) \ dy,
$$
and suppose $ \ J_{\beta}(K)  > 0. \ $  It is well known, see e.g. \cite{Compte}, \cite{Lerasle} that

$$
|{\bf E} f_n - f|  \le C_1 h^{\beta}\cdot J_{\beta}(K), \ \Var  \{f_n \} \ \le C_2 V_2(K)/nh.
$$
 \ The mean square error  for the considered statistics $ \ f_n \ $  allows the estimate

$$
Z_n(K) \stackrel{def}{=} E(f_n - f)^2 \le C_3 \left[ \ \frac{V_2(K)}{nh} + h^{2 \beta} J^2_{\beta}(K)  \ \right].
$$
\ The minimum  $ \ W = W(K) \ $ of the right - hand side of the last inequality relative the bandwidth $ \ h \ $ is following

 $$
 W = W(K) := \min_{h > 0} Z_n(K) \asymp C_4 \ n^{-2 \beta/(2 \beta + 1)} \ J_{\beta}(K) \cdot \left[ \ V_2(K)  \ \right]^{2 \beta}.
 $$
 \ We get to the following extremal problem relative the kernel $ \ K(\cdot) \ $ under formulated before limitations

\begin{equation} \label{Phi beta}
\Phi(K) \stackrel{def}{=} J_{\beta}(K) \cdot \left[ \ V_2(K)  \ \right]^{2 \beta} \to \min_K.
\end{equation}

\vspace{4mm}

 \ Let us apply the famous calculus of variations, see e.g. \cite{Courant Hilbert}, p.169; \cite{Gelfand Fomin},
 chapter 2, section 2.2. Namely, introduce the perturbed  kernel

$$
K_{\delta}(y) := K_0(y) + \delta \ g(y),
$$
where $ \ K_0 \ $ is the optimal kernel, $ \ \delta \ $ is "small" constant, e.g.
 $ \ -0.5 \le \delta \le 0.5, \ g = g(y), \ |y| \le  \theta = \theta_0 \ $ is suitable perturbation function. Of course,

\begin{equation} \label{perturb}
\int_{-\theta}^{\theta} g(y) \ dy = 0,
\end{equation}
"centering" condition. \par

 \ We obtain after some calculations

$$
\Phi \left(K_{\delta} \right) =  \Phi \left(K_{0} \right) +
$$

$$
C_2 \delta \int_{-\theta}^{\theta} \left\{ \ C_3 K_0(y) + C_4 |y|^{\beta} \ \right\} \ g(y) \ dy + 0(\delta^2), \ \delta \to 0.
$$

 \ Therefore, for some  finite constants $ \ C_3, C_4 \ $ and for arbitrary perturbation "centered" function $ \ g = g(y) \ $
 (\ref{perturb})

\begin{equation} \label{int o}
 \int_{-\theta}^{\theta} \left\{ \ C_3 K_0(y) + C_4 |y|^{\beta} \ \right\} \ g(y) \ dy = 0.
\end{equation}

 \ It follows  from  (\ref{int o}) taking into account the "centering" condition (\ref{perturb}) that

\begin{equation} \label{form of K}
 C_3 K_0(y) + C_4 |y|^{\beta} = C_5.
\end{equation}

 \ We conclude that the considered in this section optimization problem quite coincides with considered in the third section! \par

\vspace{3mm}

 \ Thus, the optimal kernel $ \ K_0 \ $  in this statement problem is described completely in the theorem 3.1. \par

\vspace{5mm}

\section{Estimation of the derivatives for density.}

\vspace{5mm}

  \hspace{3mm} It is interest in our opinion to find the optimal kernels for the problem of {\it derivative
 density estimations,} in the spirit, for example, \cite{Hansen}, pp. 12 - 16; where are described also some applications.\par

\vspace{3mm}

 \ Denote

$$
f^{(r)}(x) := \frac{d^r f(x)}{dx^r}, \ r = 1,2,\ldots;
$$
and we want to build the kernel estimation $ \ f_n^{(r)}(x) \ $ for the derivative $ \ f^{(r)}. \ $  \par

  \ We suppose that  for some $ \ m = 1,2,\ldots \ $ the density function $ \ f(\cdot) \ $  is $ \ r + 2m \ $  times continuous bounded
differentiable:

$$
\sup_{x \in R} | \ f^{(r + 2m)}(x) \ | < \infty.
$$

 \ As for the kernel $ \ K(\cdot); \ $ we assume in addition to the foregoing restrictions
  that it belongs to the Sobolev's space $ \ W_{2,r}(-\theta,\theta): \ $

\begin{equation} \label{Sobolev}
V_{r,2}(K) = \int_{\theta}^{\theta} \left[ \  K^{(r)}(y)  \ \right]^2 \ dy < \infty,
\end{equation}
 and as ordinary

$$
|y| \ge \theta \ \Rightarrow K(y) = 0; \ \int_{-\theta}^{\theta} K(y) \ dy = 1;
$$

$$
\forall s =1,2,\ldots, 2 m - 1 \ \Rightarrow \int_{-\theta}^{\theta} y^s \ K(y) \ dy = 0;
$$

$$
\int_{-\theta}^{\theta} y^{2m} \ K(y) \ dy = 1; \ K(-y) = K(y).
$$

\vspace{3mm}

 \ The kernel estimate for the derivative $ \ f^{(r)}(x) \ $ has a form

\begin{equation} \label{ker der}
f_n^{(r)}(x) = \frac{1}{n \ h^{1 + r}} \ \sum_{i=1}^n K^{(r)} \left( \ \frac{\xi_i - x}{h}  \ \right),
\end{equation}
where as before $ \ n \to \infty \ \Rightarrow h = h(n) \to 0, \ n h^{1 + r} \to \infty.  \ $\par

 \ It is known, see \cite{Hansen}, p.11 - 15 that the bias as $ \ n to \infty \ $ of $ \ f_n^{(r)}(x) \ $ under
 our condition has a form

$$
{\bf E} f_n^{(r)}(x) - f(x) \sim  C_1(f) \ h^{2m} \ \int_{-\theta}^{\theta} y^{2m} \ K(y) \ dy,
$$
and the variance may be evaluated as follows

$$
\Var \left[ \ f_n^{(r)}(x)  \ \right] \sim C_2(f) \ \frac{1}{n h^{1 + r}} \ \left[ \ V_{r,2}(K)  \ \right].
$$

\vspace{3mm}

 \ We get as before to the following extremal problem under our conditions

\begin{equation} \label{der opt}
\Phi(K) := \int_{-\theta}^{\theta} \left[ \ K^{(r)}(y)  \ \right]^2 \ dy \to \min.
\end{equation}

\vspace{3mm}

 \ The solution of this problem is quite alike to one in the second section. \par

\vspace{4mm}

 \ {\bf Theorem 5.1.} The optimal kernel $ \ K_r(y) \ $ for the considered in this section is unique and has a form

\begin{equation} \label{der ker optim}
K_r(y) = \frac{1}{2 \theta} - \frac{1}{2 \theta} P_{2m + 2r} \left( \ \frac{y}{\theta} \ \right), \ |y| \le \theta,
\end{equation}

$ K_r(y) = 0, \ |y| > \theta,  \ $ where

\begin{equation} \label{theta opt der}
\theta = [1 - \mu(2m + 2r)]^{-1/(2m + 2 r)}.
\end{equation}

 \vspace{5mm}

\section{Concluding remarks.}

 \vspace{5mm}

 \hspace{3mm} {\bf A.} The multivariate version of the kernel, in particular, optimal one, has a factorizable form

$$
V(x_1,x_2,\ldots,x_d) = \prod_{j=1}^d  K(x_j), \ d = 2,3,\ldots,
$$
see \cite{Bertin1}, \cite{Bertin2}, \cite{Devroe}, \cite{Epanechnikov}, \cite{Ibragimov Khasminskii 1},  \cite{Ibragimov Khasminskii 2},
\cite{Lerasle} etc. \par

\vspace{3mm}

 \ {\bf B.} At the same optimization problem for density measurement appears for the  so - called {\it recursive} Wolverton - Wagner's density
 estimation,  see \cite{Devroe}, \cite{Khardani}, \cite{Nadaraya}, \cite{Slaoui}, \cite{Wegman}, \cite{Wolverton Wagner} and so one. \par

\vspace{3mm}

 \ {\bf C.}  Offered here method may be generalized perhaps on the so called regression problem, i.e. when

$$
\eta_i = f(x_i) + \epsilon_i, \ i = 1,2,\ldots,n;
$$
see \cite{Gasser}, \cite{Hardle}, p.64, Theorem 3.1.\par

\vspace{3mm}

 \ {\bf D.} The case when the r.v. $ \ \xi_i \ $ are positive, may be reduced to the considered here by a transform

$$
\eta_i := \ln \xi_i,
$$
therefore

$$
f_{\eta}(x) = e^x \ f_{\xi} \left( e^x \right), \ x \in (-\infty, \infty),
$$
see \cite{Charpentier}. The case when $ \ \xi_i \in (a,b) \ $ may be considered quite analogously. \par

\vspace{6mm}

\vspace{0.5cm} \emph{Acknowledgement.} {\footnotesize The first
author has been partially supported by the Gruppo Nazionale per
l'Analisi Matematica, la Probabilit\`a e le loro Applicazioni
(GNAMPA) of the Istituto Nazionale di Alta Matematica (INdAM) and by
Universit\`a degli Studi di Napoli Parthenope through the project
\lq\lq sostegno alla Ricerca individuale\rq\rq (triennio 2015 - 2017)}.\par

\vspace{4mm}

 \ The second author is grateful to  Yousri Slaoui for sending its a very interest
 article  \cite{Khardani}. \par

\end{document}